\documentclass[12pt]{article} 

\usepackage{graphicx}  
\usepackage{amsmath}
\usepackage[margin=3cm]{geometry}
\usepackage{amsthm}
\usepackage{amsfonts}
\title{Pentagram Map by Euclidean Approach}
\author{Yusaku Mori}
\date{\today}

\begin{document}
\maketitle

\begin{abstract}
As one type of incidence theory, the geometry of pentagram map seems quite classical at first. However, this is an excellent example of such a classical idea developed into a marvellous insight by some modern approach. We introduce an alternative approach based on the setting of the ordinary Euclidean geometry.
\end{abstract}

\section{Introduction}
Given a convex $n$-gon $V=(v_1,v_2,...,v_n)$, the pentagram map is a transformation defined by joining every other vertex of $V$. We denote this transformation by $T$, acting on $V$ to produce another convex $n$-gon $U$. Figure 1 provides the visualization of this transformation. Each vertex $u_i$ of $U$ is determined as the intersection of two shortest diagonals $(v_{i-1}v_{i+1})$ and $(v_iv_{i+2})$ of $V$.\\
\\
\begin{center}
 \includegraphics[scale=0.6]{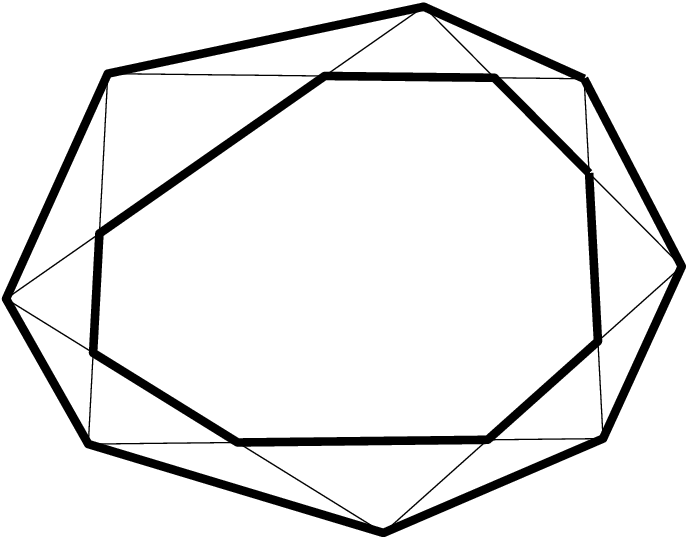}
\end{center}
\begin{center}
Figure 1.\\
The visualisation of pentagram map
\end{center}
\pagebreak
At first glance, the pentagram map itself appears to be an elementary geometry. However, in 1992, R. Schwartz published a paper [1] in which he initiated the modern approach toward the pentagram map. One of his first considerations is about the invariant of the map. He also explores the limiting flow of the pentagram map. He presents an evolution equation for this flow. Also in [2], it was shown that the pentagram map is a discrete integrable system. Due to this discovery, the study of pentagram map has become more popular nowadays.\\
\\
In section 2, we show an alternative proof of the exponential convergence to the limit point of pentagram map. In section 3, we prove the formula of the limiting flow of pentagram map, presenting that the evolution equation Schwartz showed in [1] has an extra correction term. 
\section{The exponential convergence to the limit point }
The invariant of pentagram map on the general convex $n$-gon was first discussed in [1]. There are several proofs for this invariant. See [1] and [3], for example. We are now interested in providing an alternative proof for the formula of such an invariant.\\
\\
Define
$$
f(V)=\prod_{i=1}^nX_{v_i}(V)=\prod_{i=1}^n\frac{|v_{i-1}-u_{i-1}||u_i-v_{i+1}|}{|v_{i-1}-u_i||u_{i-1}-v_{i+1}|},
$$
where $|v-u|$, $v,u\in \mathbb{R}^2$, is the signed distance between 2-dimensional vectors $u$ and $v$. Then we have the following.\\
\\
\textbf{Theorem 2.1.(Schwartz)} $f(V)=T(f(V))$. That is, the product $f$ is invariant under the pentagram map.
\pagebreak
\subsection{The Lemmata and the proof of the theorem of invariant}

To present our alternative approach, we  start with introducing some notations.\\
\\
\textbf{Definition 2.2.}
Let $V=(v_1,v_2,..,v_n)$ be a convex $n$-gon for $v_i\in \mathbb{R}^2$. If we write, 
$[u,v]:=\det
\begin{pmatrix}
u_1&v_1\\
u_2&v_2
\end{pmatrix}
$, for $u,v\in \mathbb{R}^2$, then define\\
$$
A_i=\frac{[v_{i-1}-v_{i+1},v_{i+1}-v_{i+2}]}{[v_{i-1}-v_{i+1},v_i-v_{i+2}]},
B_i=\frac{[v_{i-1}-v_i,v_i-v_{i+1}]}{[v_{i-1}-v_{i+1},v_i-v_{i+2}]},
$$
$$
C_i=\frac{[v_i-v_{i+1},v_{i+1}-v_{i+2}]}{[v_{i-1}-v_{i+1},v_i-v_{i+2}]},
D_i=
\frac{[v_{i-1}-v_i,v_i-v_{i+2}]}{[v_{i-1}-v_{i+1},v_i-v_{i+2}]}
$$
We assume that the vertices of $V$ are in general position, i.e. no three points lie on the same line, so that the determinants are non-zero. With these notations, we first present our Euclidean version of the formula for pentagram map.\\
\\
\textbf{Lemma 2.3.} Let $V=(v_1,...,v_n)$ be a convex $n$-gon, $(u_1,u_2,...,u_n):=T(v_1.v_2,...v_n)$ be the pentagram map, and $A_i,B_i,C_i,D_i$ be as above. Then, 
$$
u_i=A_iv_i+B_iv_{i+2}=C_iv_{i-1}+D_iv_{i+1}.
$$
\begin{proof}
We prove the first equality. The proof for the second one is similar. By noting that $\frac{1}{2}[u-v,v-w]$ is the signed area of $\triangle(uvw)$, for $u,v,w\in \mathbb{R}^2$, we get
\begin{equation}
\frac{B_i}{A_i}=\frac{[v_{i-1}-v_i,v_i-v_{i+1}]}{[v_{i-1}-v_{i+1},v_{i+1}-v_{i+2}]}=\frac{(1/2)(a+b)|u_i-v_i|}{(1/2)(a+b)|v_{i+2}-u_i|}=\frac{|u_i-v_i|}{|v_{i+2}-u_i|},
\end{equation}
where $a$ is the length of the line which starts at $v_{i-1}$ and extends perpendicularly to the line $(v_iv_{i+2})$, and $b$ is the length of the line starting at $v_{i+1}$ to $(v_iv_{i+2})$ in the same manner. See Figure 2. By the definition of $u_i$, we have
$$
\begin{aligned}
u_i&=v_i+\frac{|u_i-v_i|}{|v_{i+2}-v_i|}(v_{i+2}-v_i)=\frac{|v_{i+2}-u_i|}{|v_{i+2}-u_i|+|u_i-v_i|}v_i+\frac{|u_i-v_i|}{|v_{i+2}-u_i|+|u_i-v_i|}v_{i+2}\\
&=\frac{1}{1+\frac{|u_i-v_i|}{|v_{i+2}-u_i|}}v_i+\frac{1}{1+\frac{|v_{i+2}-u_i|}{|u_i-v_i|}}v_{i+2}.
\end{aligned}
$$
Now by (1), we conclude
$$
\frac{1}{1+\frac{|u_i-v_i|}{|v_{i+2}-u_i|}}v_i+\frac{1}{1+\frac{|v_{i+2}-u_i|}{|u_i-v_i|}}v_{i+2}=\frac{1}{1+\frac{B_i}{A_i}}v_i+\frac{1}{1+\frac{A_i}{B_i}}v_{i+2}=A_iv_i+B_iv_{i+2}.
$$
\end{proof}
\pagebreak
\begin{center}
\includegraphics[scale=0.6]{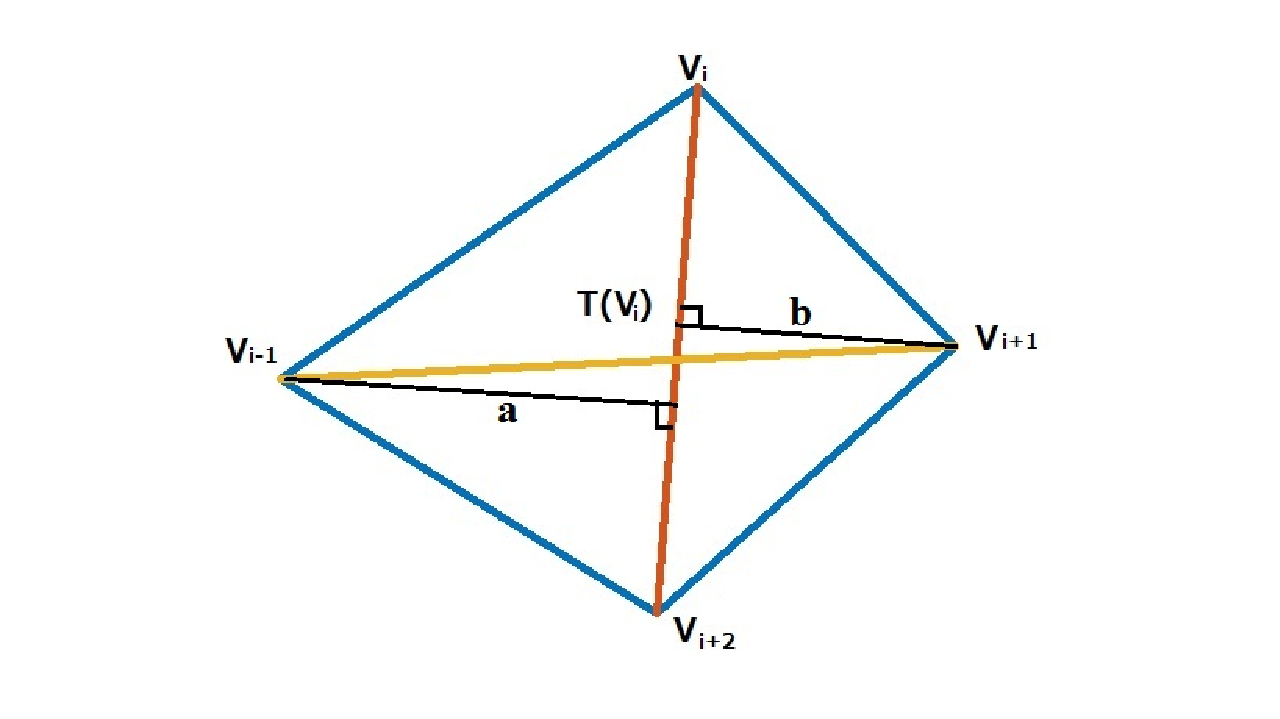}
\end{center}
\begin{center}
Figure 2. \\
Geometry of Lemma 2.3.
\end{center}
\bigskip
With this lemma, we have the following formula of how pentagram map acts on the ratio of $A_i$ and $C_i$, and $B_i$ and $D_i$.\\
\\
\textbf{Lemma 2.4.} Let $T$ be the pentagram map and $A_i$, $B_i$, $C_i$, and $D_i$ be as above. Then we have
$$
T\Big(\displaystyle\frac{C_i}{A_i}\Big)=\frac{C_i}{A_{i-1}}, \ T\Big(\displaystyle\frac{B_i}{D_i}\Big)=\frac{B_{i+1}}{D_{i+2}}.
$$
\begin{proof}
We compute the first equality. The second one is similar.
$$
\begin{aligned}
T\Big(\displaystyle\frac{C_i}{A_i}\Big)&=\frac{[u_i-u_{i+1},u_{i+1}-u_{i+2}]}{[u_{i-1}-u_{i+1},u_{i+1}-u_{i+2}]}\\
&=\frac{[C_i(v_{i-1}-v_{i+1})+B_{i+1}(v_{i+1}-v_{i+3}),(D_{i+2}-B_{i+1})(v_{i+1}-v_{i+3})]}{[A_{i-1}(v_{i-1}-v_{i+1})+B_{i+1}(v_{i+1}-v_{i+3}),(D_{i+2}-B_{i+1})(v_{i+1}-v_{i+3})]}\\
&=\frac{C_i}{A_{i-1}}
\end{aligned}
$$
\end{proof}
We are now ready to prove Theorem 2.1.
\begin{proof}[Proof of Theorem 2.1.]
We rewrite the product $f$ by applying (1) to the ratios in $f$. Then we have
$$
f(V)=\prod_{i=1}^n\frac{|v_{i-1}-u_{i-1}||u_i-v_{i+1}|}{|u_{i-1}-v_{i+1}||v_{i-1}-u_i|}= \prod_{i=1}^n\frac{B_{i-1}C_i}{A_{i-1}D_i}.
$$
Now we use Lemma 2.4. and the periodicity of the indices to get
$$
\prod_{i=1}^n\frac{B_{i-1}C_i}{A_{i-1}D_i}=\prod_{i=1}^nT\Big(\frac{C_i}{A_i}\Big)T\Big(\frac{B_{i-2}}{D_{i-2}}\Big)=T\Big(\prod_{i=1}^n\frac{B_{i-1}C_i}{A_{i-1}D_i}\Big)=T(f(V)))
$$
\end{proof}
\section{The limiting flow of pentagram map}
\subsection{The main statement on the evolution equation}
In [1], Schwartz considers the limiting flow of pentagram map. He gave a discussion on the evolution equation in [1] without any proof. Here we apply Lemma 2.3 in the previous section to give more solid argument.\\
\\
Suppose that we have a smooth parametrized curve  $\gamma: \mathbb{R} \to \mathbb{R}^2$ that is periodic. Suppose also that an $n$-gon $V=(v_1,...,v_n)$ can be expressed by $\gamma$ such that $v_{i}=\gamma(i/n)$. As Schwartz indicates, if the pentagram map converges to a flow $\gamma_t$, then we have
\begin{equation}
\frac{d\gamma_t(i/n)}{dt}\Big|_0=\lim_{n\to \infty}n^2(T^2(\gamma((i-1)/n))-\gamma(i/n)).
\end{equation}
Schwartz claims that this limit has the following evolution equation
$$
\lim_{n\to \infty}n^2(T^2(\gamma((i-1)/n))-\gamma(i/n))=\gamma''(i/n)-\frac{2}{3}W(i/n)\gamma'(i/n).
$$
where
$$
W=\frac{\det(\gamma',\gamma''')}{\det(\gamma',\gamma'')}.
$$
He first claims that $T^2(\gamma((i-1)/n))$ can be approximated by the intersection $p_i$ of the lines $(v_{i-2}v_{i+1})$ and $(v_{i-1}v_{i+2})$. However, we noticed that this is not a good approximation, and the evolution equation is slightly off by a constant, which does not vanish as $n \to \infty$. In fact, we introduce the following equation as our final result.\\
\\
\textbf{Theorem 3.1.}
$$
\frac{d\gamma_t(i/n)}{dt}\Big|_0=\frac{3}{4}\gamma''(i/n)-\frac{1}{8}W(i/n)\gamma'(i/n),
$$
\begin{flushleft}
With this equation, we arrive at the precise expression of the relation between $T^2(v_{i-1})$ and the intersection $p_i$ at the end of this section. In the following, we write $\gamma=\gamma(i/n)$, and $W=W(i/n)$.
\end{flushleft}
In order to prove theorem 3.1, we have to explore the behaviour of $T^2$. We achieve this by first showing that $A_i$, $B_i$, $C_i$, and $D_i$ do not change much by shifting the index by a small number, and later we show that the map $T$ does not change $A_i$, $B_i$, $C_i$, and $D_i$ very much.
\subsection{Lemmata}
The following lemma is the first step of the procedure.
\begin{flushleft}
\textbf{Lemma 3.2.} For bounded $k\in \mathbb{R}$, we have
\end{flushleft}
\begin{enumerate}

\item
$B_{i+k}=\displaystyle\frac{1}{4}-\frac{1}{8n}W+O\Big(\frac{1}{n^2}\Big),$
\item
$C_{i+k}=\displaystyle\frac{1}{4}-\frac{1}{16n}W+O\Big(\frac{1}{n^2}\Big) \text{ as }n\to \infty$.
\end{enumerate}
As we already saw, $D_{i+k}=1-C_{i+k}$ and $A_{i+k}=1-B_{i+k}$.
\begin{proof}
We show the calculation for $B_{i+k}$. The rest of the cases are similar.
We have
\begin{equation}
B_{i+k}=\displaystyle\frac{[v_{i+k-1}-v_{i+k},v_{i+k}-v_{i+k+1}]}{[v_{i+k-1}-v_{i+k+1},v_{i+k}-v_{i+k+2}]}.
\end{equation}
Now we express the following vectors with derivatives of $\gamma$, using Taylor expansions,
$$
v_{i+k-1}-v_{i+k}=\displaystyle-\frac{\gamma'}{n}-\frac{2k-1}{2n^2}\gamma''-\frac{3k^2-3k+1}{6n^3}\gamma'''+O\Big(\frac{1}{n^4}\Big),
$$
$$
v_{i+k}-v_{i+k+1}=-\frac{\gamma'}{n}-\frac{2k+1}{2n^2}\gamma''-\frac{3k^2+3k+1}{6n^3}\gamma'''+O\Big(\frac{1}{n^4}\Big),
$$
$$
v_{i+k-1}-v_{i+k+1}=-\frac{2\gamma'}{n}-\frac{4k}{2n^2}\gamma''-\frac{6k^2+2}{6n^3}\gamma'''+O\Big(\frac{1}{n^4}\Big), 
$$
$$
v_{i+k}-v_{i+k+2}=-\frac{2\gamma'}{n}-\frac{4k+4}{2n^2}\gamma''-\frac{6k^2+12k+8}{6n^3}\gamma'''+O\Big(\frac{1}{n^4}\Big).
$$
We plug these vectors into (3) to obtain
$$
B_{i+k}=\frac{\frac{1}{4}+\frac{k}{4n}
W}{1+\frac{2k+1}{2n}W}+O\Big(\frac{1}{n^2}\Big)=\frac{1}{4}-\frac{1}{8n}W+O\Big(\frac{1}{n^2}\Big).
$$
\end{proof}
As an important corollary, we have the following.\\
\\
\textbf{Corollary 3.3.} For bounded $k\in \mathbb{R}$, we have
\begin{enumerate}
\item $A_{i+k}=A_i+O\Big(\displaystyle\frac{1}{n^2}\Big),$
\item $B_{i+k}=B_i+O\Big(\displaystyle\frac{1}{n^2}\Big),$

\item$C_{i+k}=C_i+O\Big(\displaystyle\frac{1}{n^2}\Big),$
\item $D_{i+k}=D_i+O\Big(\displaystyle\frac{1}{n^2}\Big), \text{ as }n\to \infty$.
\end{enumerate}
\begin{proof}
This immediately follows from Lemma 2.2.
\end{proof}
\bigskip
Now we show an estimate for the map $T$.\\
\\
\textbf{Lemma 3.4.} We have
\begin{enumerate}
\item $T(A_i)=\displaystyle A_i+O\Big(\frac{1}{n^2}\Big),$
\item $T(B_i)=\displaystyle B_i+O\Big(\frac{1}{n^2}\Big),$
\item $T(C_i)=\displaystyle C_i+O\Big(\frac{1}{n^2}\Big),$
\item $T(D_i)=\displaystyle D_i+O\Big(\frac{1}{n^2}\Big), \text{ as }n\to \infty.
$
\end{enumerate}
\begin{proof}
We show the computation for $B_i$. The other cases are similar. We again apply Lemma 2.3. to get the following expression
$$
\begin{aligned}
T(B_i)&=\displaystyle \frac{[u_{i-1}-u_i,u_i-u_{i+1}]}{[u_{i-1}-u_{i+1},u_i-u_{i+2}]}=\frac{(D_i-B_{i-1})B_{i+1}}{A_{i-1}D_{i+2}-B_{i+1}C_i}\\
&=\frac{(D_i-B_{i-1})B_{i+1}}{A_{i-1}(D_{i+2}-B_{i+1})-B_{i+1}(C_i-A_{i-1})}=\frac{(D_i-B_{i-1})B_{i+1}}{A_{i-1}(D_{i+2}-B_{i+1})+B_{i+1}(D_i-B_{i-1})}.
\end{aligned}
$$
Then we apply the corollary 3.3 to establish
$$
T(B_i)=\frac{B_i}{A_i+B_i}+O\Big(\frac{1}{n^2}\Big)=B_i+O\Big(\frac{1}{n^2}\Big).
$$
\end{proof}
\bigskip
\subsection{Proof of the theorem on the flow}
Finally we prove Theorem 3.1.
\begin{proof}[Proof of Theorem 3.1.]
We begin by applying our Lemma 2.3.
$$
\begin{aligned}
T^2(v_{i-1})&=T(C_{i-1}v_{i-2}+D_{i-1}v_i)=T(v_i+C_{i-1}(v_{i-2}-v_i))\\
&=T(v_i)+T(C_{i-1})(T(v_{i-2})-T(v_i))\\
&=v_i+B_i(v_{i+2}-v_i)+T(C_{i-1})(v_{i-2}+B_{i-2}(v_i-v_{i-2})-v_i-B_i(v_{i+2}-v_i))\\
&=v_i+T(D_{i-1})B_i(v_{i+2}-v_i)+T(C_{i-1})A_{i-2}(v_{i-2}-v_i).
\end{aligned}
$$
Since $v_{i+2}-v_i=2\gamma'/n+2\gamma''/n^2+O(1/n^3)$ and $v_{i-2}-v_i=-2\gamma'/n+2\gamma''/n^2+O(1/n^3)$ as $n\to \infty$, applying Corollary 3.3 and Lemma 3.4 gives us
$$
T^2(v_{i-1})=v_i+D_iB_i(\frac{2\gamma'}{n}+\frac{2\gamma''}{n^2})+C_iA_i(\frac{-2\gamma'}{n}+\frac{2\gamma''}{n^2})+O\Big(\frac{1}{n^3}\Big).
$$
Applying Lemma 3.2 to $A_i$, $B_i$ ,$C_i$, and $D_i$, we conclude
$$
\begin{aligned}
n^2(T^2(v_{i-1})-v_i)&=n^2\Big(\frac{3}{4}+\frac{1}{16n}W\Big)\Big(\frac{1}{4}-\frac{1}{8n}W\Big)\Big(\frac{2\gamma'}{n}+\frac{2\gamma''}{n^2}\Big)\\
&+n^2\Big(\frac{1}{4}-\frac{1}{16n}W\Big)\Big(\frac{3}{4}+\frac{1}{8n}W\Big)\Big(\frac{-2\gamma'}{n}+\frac{2\gamma''}{n^2}\Big)+O\Big(\frac{1}{n}\Big)\\
&=\frac{3}{4}\gamma''-\frac{1}{8}W\gamma'+O\Big(\frac{1}{n}\Big).
\end{aligned}
$$
\end{proof}
\begin{flushleft}
\textbf{Remark}: We can also prove the following asymptotic formula for the intersection point $p_i$ of the lines $(v_{i-2}v_{i+1})$ and $(v_{i-1}v_{i+2})$ with a similar method as we saw in the proof of Theorem 3.1,
\end{flushleft}
\begin{equation}
n^2(p_i-r(i/n))=\displaystyle \gamma''-\frac{2}{3}W\gamma'+O\Big(\frac{1}{n}\Big), \text{ as }n\to \infty
\end{equation}
which is the evolution equation that Schwartz introduced in [1].\\
\\
With this remark, we have our final statement about the relation between the twice iterated pentagram map and the intersection $p_i$.\\
\\
\textbf{Corollary 3.5.} Let $p_i$ be the point of the intersection between $(v_{i-2}v_{i+1})$ and $(v_{i-1}v_{i+2})$. Then we have
$$
T^2(v_{i-1})= \displaystyle p_i-\frac{\gamma''}{4n^2}+\frac{13W\gamma'}{24n^2}+O\Big(\frac{1}{n^3}\Big), \text{ as }n\to \infty.
$$
\clearpage
\begin{center}
 \includegraphics[scale=0.4]{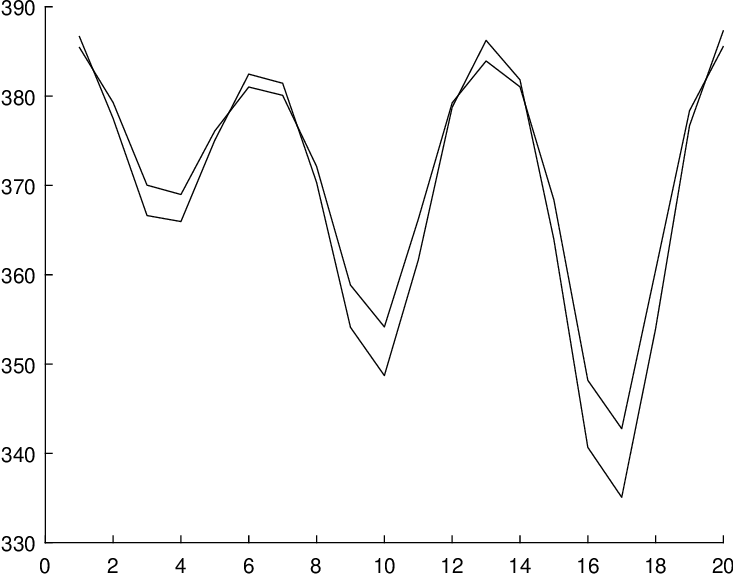}
 \includegraphics[scale=0.4]{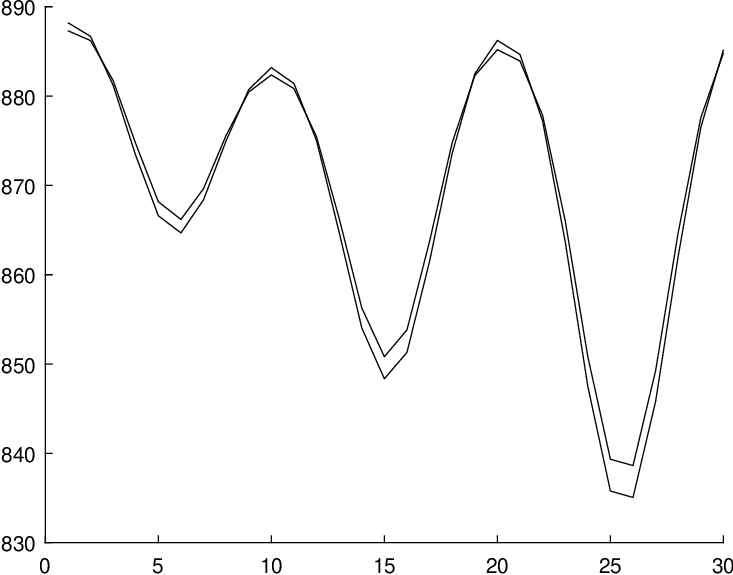}
 \includegraphics[scale=0.4]{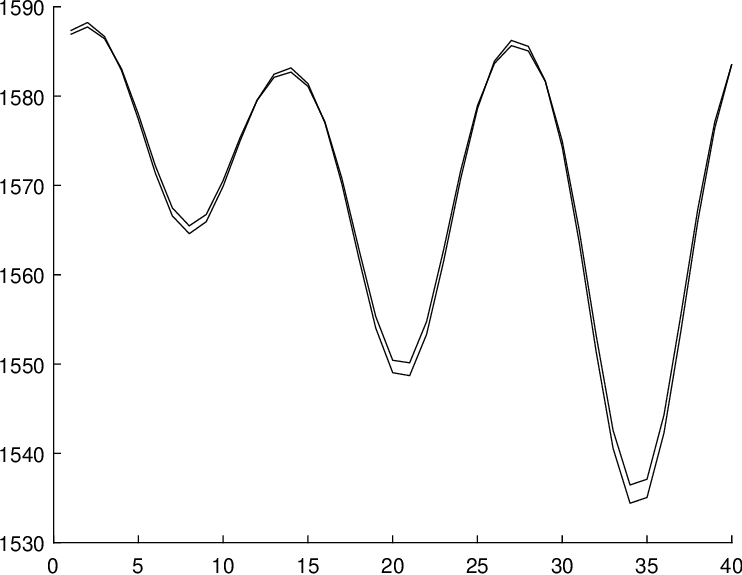}
\end{center}
\begin{center}
Figure 3
\end{center}
\begin{center}
\begin{minipage}{10cm}
Demonstration of the theorem 3.1 with $\gamma(x)=(\cos(\theta(x)), \sin(\theta(x)))$ where $\theta(x)=2\pi x+0.1\cos(2\pi x)+0.07\sin(4\pi x+\pi/3)+0.1\cos(6\pi x+\pi/5)$ for $n=20, \ 30, \text{ and } 40$.
The two curves $|n^2T^2(v_{i-1})|$ and $|n^2\gamma+(3/4)\gamma''-(1/8)W\gamma'|$ converge to the same curve as $n$ increases.
\end{minipage}
\end{center}
\medskip
\begin{center}
 \includegraphics[scale=0.4]{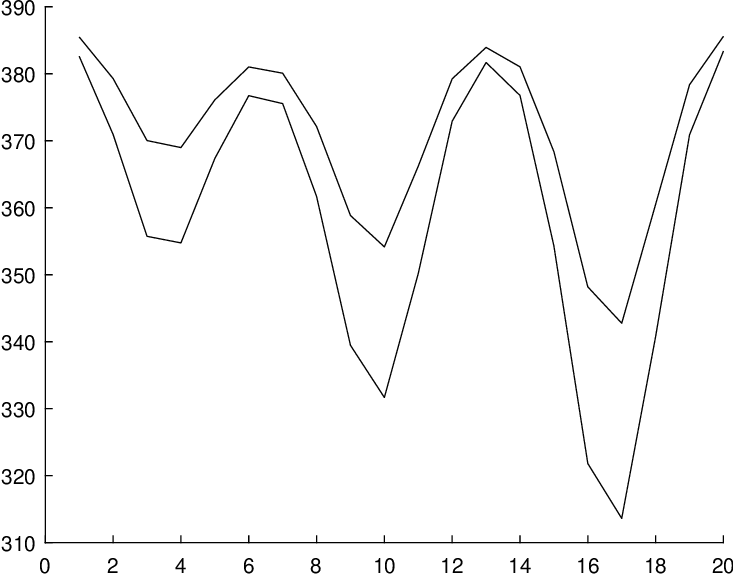}
 \includegraphics[scale=0.4]{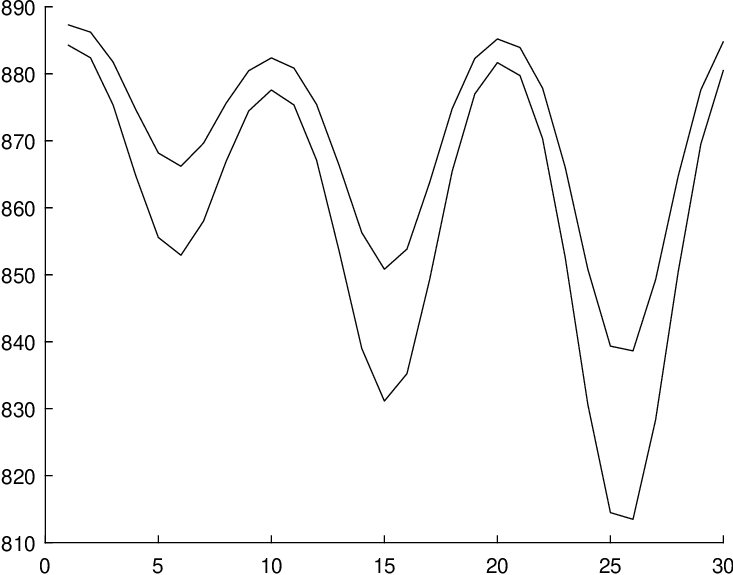}
 \includegraphics[scale=0.4]{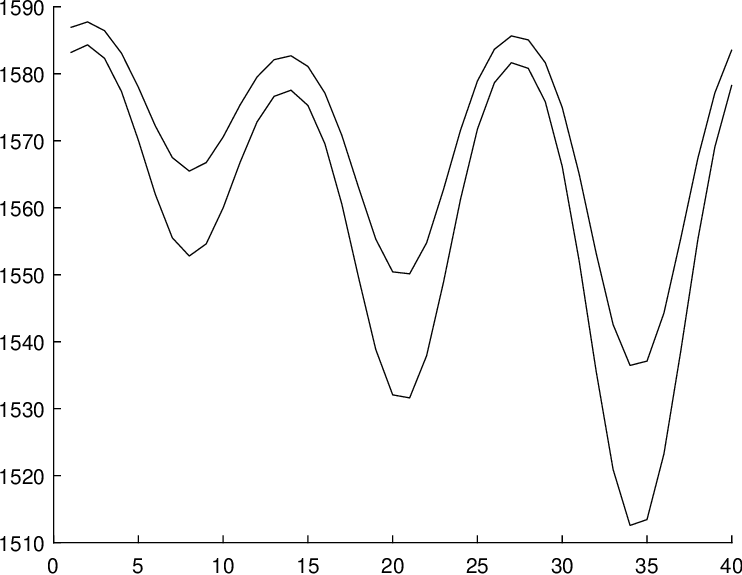}
\end{center}
\begin{center}
Figure 4
\end{center}
\begin{center}
\begin{minipage}{10cm}
The graph of two curves $|n^2T^2(v_{i-1})|$ and $|n^2\gamma+\gamma''-(2/3)W\gamma'|$ with the same $\gamma$ as Figure 3: the upper curve is the curve for $|n^2T^2(v_{i-1})|$, and the lower for $|n^2\gamma+\gamma''-(2/3)W\gamma'|$. We see that the two curves do not converge to the same curve, but keep a certain distance even if $n$ becomes large.
\end{minipage}
\end{center}

\clearpage
\pagebreak
\section{Reference}
$$
\begin{aligned}
\text{[1]}.& \text{ R. Schwartz, \textit{The pentagram map}, Experiment. Math. 1 (1992), 71-81.}\\
[2].& \text{ V. Ovsienko, R. Schwartz, and S. Tabachnikov,}\\
& \text{ \textit{The pentagram map: a discrete integrable
system}, Comm. Math. Phys. 299 (2010), 409-446.}\\
[3].& \text{ R. Schwartz, \textit{Discrete monodromy, pentagrams, and the method of condensation},} \\
& \text{ J. of Fixed Point Theory and Appl. 3(2008),379-409.}\\
[4].& \text{ R.Schwartz, \textit{Recurrence of the Pentagram Map}(2000)}.\\
[5].& \text{ M. Glick, \textit{The limit point of pentagram map}.}
\end{aligned}
$$
\end{document}